\documentclass{article}
\usepackage[usenames,dvipsnames]{color}
\usepackage{amsmath}
\usepackage{amsfonts}
\usepackage{amssymb}
\usepackage{amsthm}
\usepackage{manfnt}
 \newcommand{\Ob}{\operatorname{Ob}\,}
   \newcommand{\Mor}{\operatorname{Mor}\,}

\usepackage{graphicx}
\usepackage{amscd}
\usepackage{qtree}
\def\trd{\textcolor{red}}

\def\tbi{\textcolor{Bittersweet}}

\newcommand{\eq}[1]{\begin{equation}
                     \begin{split} 1 \end{split}
                     \end{equation}}
      
      \def\h{\hskip5ex}               
                     \newcommand{\cinf}{C^{\infty}}
\newcommand{\lon }{\longrightarrow }

\def\p{\vskip2ex}
\newcommand{\A}{{\mathcal{A}}}

\newtheorem{theorem}{Theorem}[section]

\theoremstyle{definition}
\newtheorem{definition}[theorem]{Definition}

\theoremstyle{definition}

\newtheorem{dfn}[theorem]{Definition}

\theoremstyle{remark}
\newtheorem{remark}[theorem]{Remark}

\newtheorem{rmk}[theorem]{Remark}

\def\AA{$A_\infty$}
\def\LL{$L_\infty$}
\def\BBvD{Berends, Burgers and van Dam } 
\def\BBvDs{Berends, Burgers and van Dam} 

\def\XX{\Xi}
\def\PP{\Phi}
\def\LP{\Lambda^*\Phi}
\def\DD{\mathcal{D}}

\def\pp{\phi}

\def\dd{\delta}

\def\BV{Batalin-Vilkovisky}
\def\BFV{Batalin-Fradkin-Vilkovisky\ }
\def\BBvD{Berends, Burgers and van Dam } 
\def\FLS{Fulp, Lada and Stasheff}
\def\BFLS{Barnich, Fulp, Lada and Stasheff}
\def\HZ{Hohm and Zwiebach}
\def\pp{\phi}

\def\be{\begin{equation}}
\def\ee{\end{equation}}
\def\bea{\begin{eqnarray}}
\def\eea{\end{eqnarray}}
\def\p{\vskip2ex}
\def\aa{\alpha}
\def\bb{\beta}
\def\gg{\gamma}

\newcommand{\da}{_{\alpha}}   %subscript alpha
   \newcommand{\db}{_{\beta}}    %\subscript beta
   \newcommand{\dgam}{_{\gamma}}   %subscript gamma
   \newcommand{\dab}{_{\alpha\beta}}   %subscript alpha\beta
   \newcommand{\dagam}{_{\alpha\gamma}}   %subscript alpha\gamma
   \newcommand{\dbgam}{_{\beta\gamma}}

\def\MB{{\boldsymbol \Omega} B}
\def\MY{{\boldsymbol \Omega} Y}
\def\MX{{\boldsymbol \Omega} X}

\def\LO{Loday}
\def\LOO{$Loday_\infty$}
  \def\D{\Delta}

\title{ \LL\ and \AA-structures: then and now }
\author{Jim Stasheff}
\begin{document}
 \pagestyle{myheadings}
\markboth{\LL\  \ \&\  \AA:  then and now}{\LL\  \&\  \AA: then and now\hskip10ex\today}

\maketitle
\begin{abstract}
\ Looking back over 55 years of higher homotopy structures, I 
reminisce as I recall the early days and ponder how they developed and how I now see them. From
the history of \AA-structures and later of \LL-structures and their progeny,  I hope to highlight some old results
 which seem not to have garnered the attention they deserve as well as some tantalizing new connections.  

\end{abstract}
\vskip2ex
\noindent Dedicated to the  memory of Masahiro Sugawara and John Coleman Moore
%\vfill\eject

\tableofcontents
%\vfill\eject
\section{Introduction}

\h Looking back over 55 years of higher homotopy structures, I 
reminisce as I recall the early days and ponder how they developed and how I now see them. From
the history of \AA-structures and later of \LL-structures and their progeny, I emphasize my \emph{homotopy} perspective on
  how  they morphed and intertwined in homotopy theory 
with applications to 
geometry and physics.  A recurring theme is the relation of higher algebraic structures with higher topological, geometric or physical structures \cite{mss}.
There also important higher algebraic structures without homotopy, where `higher' here means generalizations Lie brackets to $n$-ary brackets 
for $n>2$.  I will touch on these only briefly (see Section \ref{n-ary}); a very thorough survey is provided by de Azcarraga and Izquierdo \cite{Az&Iz:n-ary}.
\p
Since 1931 (Dirac's magnetic monopole), but especially in the last six decades, there has been
increased use of cohomological and even homotopy theoretical techniques in mathematical physics.
 It all
began with Gauss in 1833, if not sooner with Kirchof's laws. The cohomology referred to in Gauss
was that of differential forms, div, grad, curl and especially Stokes Theorem (the de Rham complex). I'll mention some of the  more `sophisticated' tools now being used.
  \p
As I tried to being this survey reasonably up to date, one thing led to another, reinforcing my earlier image not of a tree but of a spider web.
I finally had to quit pursuit before I became trapped!
My apologies if your favorite strand is not mentioned.
\p
I have included bits of history with dates which are often for the published work, not for the earlier arXiv post or samizdat.
\p
Acknowledgements: I am grateful to the editors for this opportunity, to all my coauthors and to the many who have responded  to earlier versions.
The remaining gaps, of which there are several, are mine.
%, I am especially grateful to  Fusun Akman, Yvette Kosmann-Schwarzbach, and Thedia Voronov for significant help, especially with Section \ref{derived} and to the chief operadchik Bruno Vallette.
\section{Once upon a time:  \AA-spaces, algebras, maps}
\h For me, the study of associativity began as an undergraduate at Michigan. I was privileged to have a course in classical projective geometry from the eminent relativist George Yuri  Rainich\footnote{Born Yuri Germanovich Rabinovich  $https://en.wikipedia.org/wiki/George_Yuri_Rainich$}. A later course ``The theory of invariants'' emphasized the invariants of
Maxwell's equations. Rainich included secondary invariants, preparing me well when I later encountered secondary cohomology  operations.
\subsection{\AA-spaces}
\hskip5ex
The history of \AA-structures begins, implicitly, in 1957 with the work of Masahiro Sugawara \cite{sugawara:h,sugawara:g}. 
He showed that, with a generalized notion of fibration,  the  Spanier-Whitehead condition for a space $F$ to be an $H$-space:

 \emph{The existence of a fibration with  fibre contractible in the total space}
 
\noindent  is necessary and sufficient.
Sugawara goes on to obtain similar criteria for $F$ to be a homotopy associative $H$-space or a loop space.
\p
When I was a graduate student at Princeton, John Moore suggested I look at when a primitive cohomology class $u\in H^n(X,\pi)$ of a topological group or loop space $X$ was the suspension of a class in $H^{n+1}(BX,\pi)$.  In other words, when was an H-map $X\to K(\pi,n)$ induced as the loops on a map 
$BX\to K(\pi,n+1)$. Here, for a topological group or monoid, $BX$ refers to a ``classifying space''.
\p
 In part inspired by the work of Sugawara, I attacked the problem in terms of a filtration  of the classifying bundle $EX\to BX$ by the projective spaces 
$XP(n)$. Sugawara's work was for $XP(3)$ and $XP(\infty)$. 
%(More about that filtration later.)
His criteria consisted of an infinite sequence of conditions, generalizations of homotopy associativity.  He included conditions involving homotopy inverses 
which I was able to avoid via mild restrictions.
%\footnote{For me, associativity began as an undergraduate at Michigan. I was privileged to have a course in classical projective geometry from the eminent relativist George Yuri  Rainich (born Yuri Germanovich Rabinovich ) $https://en.wikipedia.org/wiki/George_Yuri_Rainich$}.  

\begin{remark} Any H-space has a projective `plane' $XP(2)$ and homotopy associativity implies the existence of $XP(3)$. Contrast this with  classical projective geometry where the existence of a projective 3-space implies strict associativity. 
\end{remark} 
To systematize the higher homotopies, I defined:

\begin{dfn}
An $A_n$ \emph{space} $X$ consists of a space $X$
 together with a coherent set of maps 
$$m_k: K_k \times X^k\to X\ \  for \ \ 2\leq k\leq n$$
where $K_k$ is the (by now) well known $(k-2)$-dimensional associahedron.
\end{dfn}
%\vfill\eject 

My original realization had little symmetry, being based on parameterizations, following Sugawara.

\begin{figure}[ht]
    \begin{center}
        \includegraphics[height=1.2in]{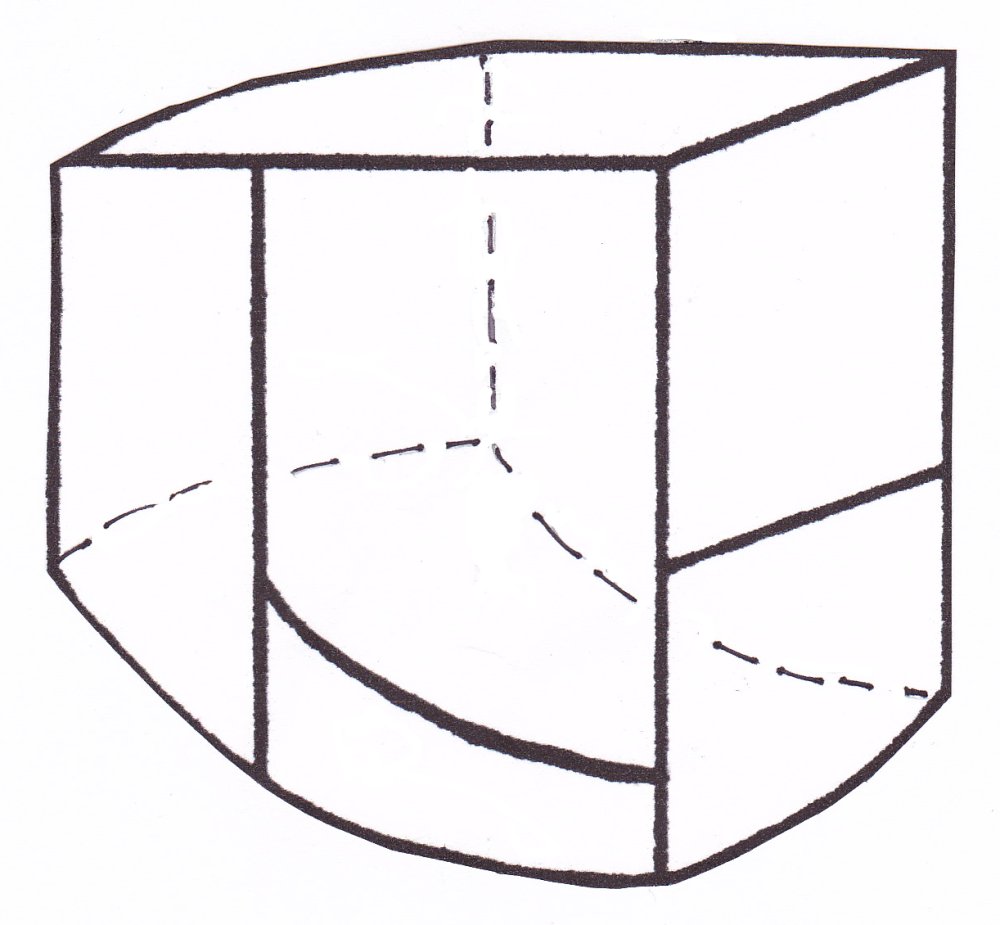}
         
    \end{center}
    \end{figure}

The \emph{``so-called''} Stasheff polytope was in fact constructed by Tamari in 1951 \cite{tamari:thesis,tamari-1962, tamari-huguet:polyedrale}, a full decade before my version.
\begin{figure}[ht]
    \begin{center}
        \includegraphics[height=2in]{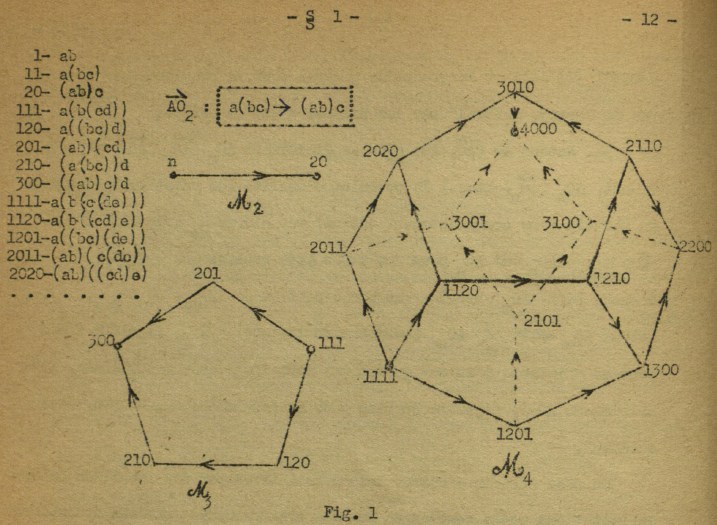}
    \end{center}
    \caption{from Tamari's thesis \cite{tamari:thesis}}
    \end{figure} 
 
 Many other realizations are now popular and have been collected by Forcey \cite{forcey:zoo} along with other relevant ``hedra''.
\p
Tamari's point of view was much different from mine, but just as inspiring for later work, primarily in combinatorics. The book \cite{tamari:book} has a wealth of offspring.   
  
\p
My multi-indexing of the cells of the associahedra was awkward, but the technology of those years made indexing by trees unavailable. 
Boardman and Vogt use
spaces of binary trees with interior edges given a length in $[0,1]$, producing a cubical subdivision of the associahedra \cite{BoVo}. Alternatively, the assocciahedron $K_n$ can be realized as the compactification of the space of $n$ distinct points in $[0,1]$ \cite{mss} Example 4.36, \cite{jds:config} Appx B. For  $n$ distinct points on the circle, the compactification
is of the form $S^1\times W_n$, where $W_n$ is called the cyclohedron (see Section \ref{quest}).
\p
Just as one construction of a classifying space $BG$ for a group $G$ is as a quotient of the disjoint union of pieces 
$\Delta^n\times G^n$ where $\Delta^n$ is the n-simplex, I constructed $BX$ as a quotient of the disjoint union of pieces $K_n \times X^n$.
%Note  this is not the same as an \AA-space with $m_k=0$ for $k > n.$  For example, for $n=3$ an associating homotopy need have no relation to the usual pentagon relation.
\p
 I was still working on Moore's problem when I went to Oxford as a Marshall Scholar. I was automatically assigned to J.H.C. Whitehead as supervisor, then transferred to Michael Barratt whom I knew from Princeton. When Barratt moved to Manchester,
Ioan James took on my supervision. Though Frank Adams was not at Oxford, he did visit and had a significant impact on my research.
Ultimately, this resulted in my theses for Oxford and Princeton. 
I needed to refer to what would be my Princeton thesis as a prequel to the Oxfod one.
Oxford had a strange rule that I could not include anything I \emph{had submitted} elsewhere for a degree, even with attribution, so I made sure to submit to Oxford before returning to the US and finishing the formalities at Princeton. For the record, the 1963 published versions \emph{Homotopy Associativity of H-spaces I and II}
\cite{jds:hahI, jds:hahII}
correspond respectively to the topology of the Princeton thesis and the homological algebra of the Oxford one.
\subsection{\AA-algebras}
\h Thinking of the cellular structure of the associahedra led to:
 \begin{dfn}
 Given a graded vector space  $A= \{A_n\}$, an \AA-algebra structure on $A$ is a coherent set of maps $m_k:A^{\otimes k}\to A$ of degree $k-2$.
 \end{dfn}
 Here \emph{coherence} refers to the relations
 \be
0 = \sum_{i+j=n+1}\sum_{k=1}\pm m_i \bigl(a_1\otimes \cdots \otimes a_{k-1}\otimes m_j( a_k\otimes\cdots\otimes a_{k+j-1}) \otimes a_{k+j}\otimes\cdots\otimes a_n\bigr)\ee
 \p
 For an ordinary dga  $(A,d,m)$, Massey constructed secondary operations now called \emph{Massey products} (Massey's $d$ was cohomological, i.e. of degree +1).
 These products generalize easily to an \AA-algebra, using $m_3$, etc.
 \p
For an ordinary associative algebra $(A,m)$, there is the \emph{bar constuction} $BA$, a differential graded coalgeba with differential determined by the multiplication $m$.
For an \AA-algebra, there is a completely analogous construction using all the $m_k$. The most effective defintion of an \AA-morphism  from $A_1$ to $A_2$ is as a morphism 
of dg coalgebras $BA_1\to BA_2$. There is an important  but subtle relation between the differentials of the \emph{bar construction spectral sequence}
and (higher) Massey products (see \cite{buijsetal:massey}).
\p
%\trd{Something about Massey and Yessam with ref to Retakh below}
\subsection{\AA-morphisms}
\h
Morphisms of \AA-spaces are considerably more subtle. Again the way was led by Sugawara \cite{sugawara:hc}, 
who presented \emph{strongly homotopy multiplicative maps} of strictly associative $H$-spaces.  
\begin{dfn}
\label{shmap}
A  \emph{strongly homotopy multiplicative map} $X\to Y$ of associative $H$-spaces consists of a coherent family of maps $I^n\times X^n\to Y$.
\end{dfn}
The analog for dg associative algebras is straightforward.
\p
When  $X$ and $Y$ are \AA-spaces, the parameter spaces for the 
higher homotopies, now known as the \emph{multiplihedra}  (often  denoted ${\mathcal J}_n$ or $J_n$),   are considerably more complicated and 
realization as convex polytopes took much longer to appear.
%\footnote{
An excellent illustrated collection of these and other related polyhera is 
given by Forcey \cite{forcey:multiplihedra}, see also \cite{Kawa}.
%YusukeKawamoto,HigherhomotopycommutativityofH-spacesandhomotopylocalizations,PacificJournal of Mathematics Vol. 231, No. 1, (2007) 103-126}.
 My early picture for $n  = 4$ was 
not really a convex polytope though  my drawing of $\cal{J}_4$ appears as a subdivided pentagonal
cylinder.  Iwase and Mimura in \cite{iwase-mimura:hha} gave  the first detailed definition of the
multiplihedra and describe their combinatorial properties.   If the range $Y$
 is strictly associative, then the multiplihedron ${\mathcal J}_n$
collapses to the associahedron $K_{n+1}$ \cite{jim:book}. On the other hand, Forcey \cite{forcey:quotient} observed that 
if the domain $X$ is strictly associative, then the multiplihedron denoted ${\mathcal J}_n$
collapses to
a \emph{composihedron} he created and denotes ${\mathcal CK}(n)$, new for $n\geq 4$.
\p
The analog for \AA-algebras is much easier to write down, convexity not being an issue.
 The study of  $A_{\infty}$ spaces and algebras continues.
There are interesting questions about the extension of
$A_n$-maps, as in \cite{hemmi:review} and about the transfer of $A_{\infty}$ structure through these
maps, as in \cite{markl:tr}.

\section{Iterated loop spaces and operads}\label{iterated}
\h
In terms of the  development of higher homotopy structures, perhaps my most important result was this characterization of spaces of the homotopy type of loop spaces in terms of an \AA-structure: 

\begin{theorem} A `nice' connected space $X$ has the homotopy type of a based loop space $\Omega Y$ for some $Y$ if and only if $X$ admits the structure of an \AA-space.
\end{theorem}

Here `nice' means of the homotopy type of a CW-complex with
a non-degenerate base point. For the standard description of $\Omega Y$ in terms of loops parameterized by $[0,1]$, the maps are generalization of the usual one for homotopy associativity.
\p
An alternative is to use the strictly associative \emph{Moore space of loops} $\MY$ with  loops parameterized by intervals $[0,r]$ for $r\geq 0$ \cite{moore:loops}\label{mooreloops}. Apparently  Moore never formally published this major contribution, but, thanks to the internet, his seminar of 1955/56 is available at 
http://faculty.tcu.edu/gfriedman/notes/aht23.pdf. The homotopy equivalences of $X$ and $\Omega Y$ and $\MY$
are indeed \AA-morphisms.
\p
From characterizing loop spaces, I went on to characterizing loop spaces  on H-spaces by constructing a multiplication on a classifying space $BX$.  This led to  issues of homotopy commutativity with formulas getting out of hand \cite{jim:hc}. Soon after,  Clark \cite{clark:hc}  investigated the subtleties of comparing ${\boldsymbol \Omega}(X\times X)$ with 
$\MX\times\MX$.
\p
 In 1967 at the University of Chicago, prompted by a visit from Frank  Adams, a seminar {\it Iterated homotopies and the bar construction} 
was organized by Adams and Mac Lane. Thanks to Rainer Vogt we have a great record of that seminar \cite{vogt:memories}. The audience was exceptional; see Vogt's listing as well as his recollection of the lectures. He writes that the seminar had a great influence on his work with Mike Boardman, as it did on several of the participants, myself included.
\p
In 1968 Boardman and Vogt \cite{BoVo:bull,BoVo} (1973), motivated by the many infinite loop spaces then of interest thanks to Bott periodicity and that seminar,
emphasized the point of view of \emph{homotopy invariant algebraic structures} to characterize such infinite loop spaces. 
\p
%\trd{integrate the follow}
A key idea was the passage from a strict algebraic structure to one that was  homotopy invariant but still of the same homotopy type, suitably interpreted. Since retiring from UNC, I have been  a long time guest at U Penn, where the language of algebraic geometry is dominant, puzzling me as to what  `derived' referred. To paraphrase Monsieur Jourdain, I was happy to  discover that I had been speaking `derived' all my life; i.e. speaking in a related homopy category (see Section \ref{derived} for an alternate meaning).
\p
Then around 1970,
Peter May came along and blew the subject wide open with his development of operads  \cite{may:geom} to handle the 
complex of homotopy symmetries to characterize iterated loop spaces of any level. Of particular importance early on for Boardman and Vogt and  for May was the little $n$-cubes operad.Somewhat later, Getzler \cite{getzler:bvand2D} introduced the little $n$-disks operad and its framed version.
Since then there has been a proliferation of operads with additional structure as well as many generalizations of the concept, \cite{loday-vallette, mss:opbook}.
%, beyond the infinite loop spaces so well handled by Boardman and Vogt. 
\p
Operads were crucial for studying an  important issue in the $\infty$-version of commutative algebras: whether to relax the commutativity up to homotopy or to keep the strict symmetry but relax the associativity or relax both.  \p
 As emphasized by Kontsevich \cite{kont:trium}, the triumvirate of \AA-, \LL-\ and $C_\infty$-algebras play a dominant role. 
 By $C_\infty$-algebra is meant what is also known as a balanced \AA-algebra, that is, 
 a strictly graded
commutative \AA-algebra defined  in terms of a coherent set of n-ary products which vanish on shuffles.  
$C_\infty$-algebras and \LL-algebras are in an adjoint relationship just as are
strict associative commutative algebras and Lie algebras. The next most prominent might be $E_\infty$-algebras, dubbed \emph{homotopy everything}. That is a bit of a misnomer, though they are very important for the study of infinite loop spaces. 
\p
As readers of this journal are well aware, there has been a proliferation of strong homotopy or $\infty$-structure
 versions of classical algebras such as Gerstenhaber, Batalin-Vilkovisky, Froebenius
% \trd{Feynman?}
and on and on  (see Sections \ref{derived} and \ref{bvalg}) as well $\infty$ `spaces'
and even $\infty$ `group/oids'.  Those occur in attempting `integration' as in integration of a Lie algebra to a Lie group \cite{henriques:int} 
regarded as a simplicial set or Kan complex. There is also \emph{derived \AA}
 for use over a ring rather than a field \cite{sagave:derived,livernet:derivedAoo}.

%T. J. Lada, Strong homotopy algebras over monads

\section{\LL-algebras}\label{Loo}
\h
%\trd{see the book as to what to say esp what generality}
In contrast to  homotopy associativity, which was considered long before the higher homotopies were recognized, algebras with just `Jacobi up to homotopy'  appeared 
%\trd{Alan and Dee - 1998}
only after the full set of higher Jacobi homotopies were incorporated in
\LL-algebras (aka strong homotopy Lie algebras or sh-Lie algebras), which in turn had waited many years to be introduced for lack of applications. 
However,  `Jacobi up to homotopy' was implicit in the many proofs of the graded Lie algebra structure of the Whitehead product
%\footnote
 since, as Massey said \cite{Haring},:
\begin{quote}This question was `in the air' among
homotopy theorists
in the early 1950's,  I don't 
  believe you can point to any one person and say that he or she raised this.
\end{quote}
 \p
  Retakh and Allday were the first to define Lie-Massey operations, both in 1977: \cite{retakh:massey} in Russian  and \cite{allday:massey}.
  A preferable name might be \emph{Massey brackets}. For Retakh, 
  they appeared as obstructions to deformations of complex singularities. 
  (Compare \cite{douady,retakh:douady} and see Section \ref{defs}.) Retakh's  $n$-homotopy multiplicative maps are a special case of  \LL-morphisms.
  Both Retakh and Allday emphasize applications to the Quillen spectral sequence and to rational Whitehead products. For higher Whitehead products, see also \cite{belchietal:hiwhitehead}.
  \p
  Clear exposition of the Massey brackets and their connection to rational homotopy theory is in chapter V of
  \cite{tanre:modeles}.
  %chapter V of Daniel Tanré?s book ?Homotopie rationelle: modeles de Chen, Quillen, Sullivan?. 
  
 \p
 Higher homotopies for Jacobi led to:
 
 \begin{definition} \cite{ls,lada-markl} An \LL-structure on a graded vector space $V$  is a collection of skew graded
symmetric linear brackets $l_n:\bigotimes ^nV \longrightarrow V$ 
$$\ell_n = [\ ,\dots,\ ]:\Lambda^n V\to V$$ for $n\geq 1$
of degree $2-n$
(for cochain complexes, and $n-2$ for chain complexes)
such that
 \be 
 \sum_{i+j=n+1}\pm l_i (l_j(v_{\sigma (1)}\otimes \dots \otimes v_{\sigma (j)}) \otimes v_{\sigma(j+1)}
\otimes \dots\otimes v_{\sigma (n)})=0
\ee
where
  $\sigma$
runs through all $(j,n-j)$ unshuffles  and for which `there exists a set of signs' (folk saying); in this case, the sign of the unshuffle in the graded sense.
\end{definition}
%(cf. \cite J).
\p
%To emphasize $\Loo$-algebras as generalization of Lie algebras, it is helpful to write
%\be
%l_n(v_1,\cdots,\v_n) =: [v_1,\cdots,\v_n]
%\ee
Equivalently, an \emph{\LL-algebra} is a graded vector space $L=\{L_i\}$ with a coderivation differential of degree $\pm 1$ 
on the graded symmetric coalgebra $C(L)$ on the shift $sL$.
For an ordinary Lie algebra, this is the classical Chevalley-Eilenberg {\em chain} complex.
\p
The map $l_1$ is a differential: $(l_1)^2=0$ and $l_2$ may be compared to  an ordinary
(graded)
Lie bracket $[\ , \ ].$ 
When  $l_1$ or $l_3=0$, the definition yields the usual
(graded) Jacobi identity.  In general, $l_3$ is a homotopy between the Jacobi
expression and $0$ while
the other $l_n$'s are  known as \emph{higher homotopies} or \emph{higher brackets}.
If $l_3=0$, there still may be non-trivial higher $l_n$, e.g. on the homology of a dg Lie algebra (see Section\ref{qht}).
%\footnote{
Physicists like to say products,
 but we have consistently used \emph{`brackets'}.
There are  alternative notations:
$$  l_n(v_1,\cdots,v_n) =   [v_1,\cdots,v_n] \  \text(math)\ \ =   [v_1\cdots v_n] \ \text(physics).$$
%In BBvD, the map $l_1$ is denoted $\partial$ with application to symmetric tensors in mind.
%A shorthand for indicating equal entries $v=v_1= \cdots
\p
%For an ordinary Lie algebra, this is the classical Chevalley-Eilenberg {\em chain} complex.

\begin{remark} Note the ambiguity as to the degree $\pm 1$ of $d$ in defining an \LL-algebra.
The binary operation is always of degree 0; sometimes the `manifest'  grading in examples is $ \bold{not}$ the right one; see examples below. The shift of the bracket now has the same degree as the shift of $\ell_1$.
Notice also this bracket extends to an action of the degree 0 piece on the piece of degree 1 (or -1 respectively), as for a module over an algebra (contrast
BBvD structures in Section \ref{BBvD}).
\end{remark}
\p
 Belatedly,  Sullivan's  models for rational homotopy types were recognized as  being of the form  $C(L)$ on the shift $sL$ of the \LL-algebra of rational homotopy groups \cite{sullivan:inf}. The 2-bracket corresponds to the Whitehead or Samelson product  and higher order brackets to higher order Whitehead products \cite{belchietal:hiwhitehead}.
 
\subsection{\LL structures  in physics}\label{LLphysics}
\h
  Giovanni Felder was all too briefly my colleague at UNC. Later he joined the
  \LL  club, remarking
  
\emph{  ``the $\infty$-virus had a long incubation time and the outbreak came
after I left the infection zone.''}

Since then, it has expanded to epidemic proportions in the field theoretic physics community.

\p
In 1982, \LL-algebras appeared in disguise in gravitational physics in work of D'Auria and Fr\'e. Unfortunately they referred to their algebras as free differential algebras; to be precise, their FDA is a dgca (free as a gca ignoring the differential).
% as   in Sullivan's models of rational homotopy types (1977) , in which $L_\infty$ structure was implicit though not named/recognized as such.
\p  
Around 1984, Gerett Burgers visited Henk Van Dam at UNC and we discussed parts of his thesis;
again it looked like an $L_\infty$ structure was lurking there. Later this was confirmed by \FLS. The essential idea in 
\BBvD \cite{burgers:diss,BBvD:probs} was a novel attack on  particles of higher spin by letting the gauge parameters act in a 
\emph{field dependent} way.
 %This led \FLS \ \cite{fls} and \BFLS \ \cite{bfls:gauge} \  to find an $\Loo$-algebra implicit in their formalism. 
 \p
 In 1987, the formulas of the BRST operator in the construction of \BFV for constrained Hamiltonian systems \cite{BF,FF,FV} could be recognized as corresponding to an \LL-algebra
 \cite{jds:shrep,jds:ghost,jds:hrcpa}  as did the \BV\  corresponding Lagrangian formulas   \cite{jds:moskva}.
 \p
In 1989, the $L_\infty$ structure of CSFT (Closed String Field Theory) was first identified when  Zwiebach fortuitously gave a talk in Chapel Hill at  the last GUT (Grand Unification Theory) Workshop \cite{z:csft,z:book}.
\p
By 1993, it seemed appropriate to provide an \emph{Introduction to sh Lie algebras for physicists}  \cite{ls}.
\p
In 1998, Roytenberg and Weinstein,
building on Roytenberg's thesis,  showed that Courant algebroids give
rise to (small) \LL-algebras (see Section \ref{def:quasi-algebroid} and \cite{zambon:hicourant}).   
\p
 \LL-algebras are %becoming increasingly 
 continuing to be useful (and popular !) in physics\footnote{See https://ncatlab.org/nlab/print/L-infinity+algebras+in+physics for an extensive, 
 annotated chronological list - thanks to Urs Schrieber- but beware the linguistic problems.} in two ways:
\begin{itemize}
\item Solution of a physical problem leads to a structure which later is recognized as that of an \LL-algebra.
\item Solution of a physical problem is attacked using knowledge of \LL-algebras.
\end{itemize}
 There are some famous `no go' theorems that rule out certain physical  models, e.g. for  higher spin particles (see Section \ref{BBvD}).
What are ruled out are only models in terms of representations of strict Lie algebras (compare Section \ref{ruths}).

  \subsection{`Small'  \LL-algebras}
\h
CSFT requires the full panoply of higher brackets of all orders, but many other examples of \LL-algebras with at most 3 pieces in the grading
have appeared in physics. More generally, the name \emph{Lie $n$-algebra} uses the $n$-categorical language and refers to an \LL-algebra $L$ with $L_n$ = 0 for 
$n<0$ and for $n> n-1$ or for $n>0$ and $n> -n+1$. Some authors refer to these as  `truncated' \LL-algebras.
\p

Zwiebach and other physicists had asked about small examples of \LL-algebras, (physicists' `toy' models) leading to work of Tom Lada and his student Marilyn Daily. She classified all 3-dimensional \LL-algebras: 2-graded with one 1-dimensional component and one 2-dimensional component; 3-graded where each component is 1-dimensional \cite{daily}. 
%Together they  showed that the 2-graded examples fit into the context of BBvD theory. 
\p

\dbend \ {\bf Warning!!} There is a real problem of nomenclature.
\p
There are also  notions of $n$-Lie algebra which have a $k$-ary bracket \emph{only} for $k=n$ (see Section \ref{n-ary}).
% Recall the $A_n$-algebras mentioned above. In the Lie context,   $L(m)$ denotes a structure on a graded vector space $X$ with higher brackets 
%$\ell_k = [\ ,\dots,\ ]$ given only for $k\leq m$ but satisfying the corresponding subset of the higher Jacobi identities. 
%Occasionally the term \emph{$n$-ary Lie} is used, as Filippov did originally! which would be much less confusing.
\p

\def\h{\hskip5ex}

In 2017, \HZ\  \cite{hz} introduced a wide class of field theories they call \emph{standard-form}, which included, by assumption,  a (usually small) $L_\infty$-structure;
again field dependence was crucial. They went beyond \BBvD by including
 a space of `field equations' in their   \LL-algebras. Together with Andreas Deser, Irina Kogan and Tom Lada, we are adding
 such field equations to the \BBvD approach to show that again field dependence implies the \LL structure.
 \p
 Following \HZ,  Blumenhagen, Brunner,  Kupriyanov and L\"{u}st \cite{bbkl}
attacked the existence of an appropriate $L_\infty$ structure with given initial terms
 by what they call a bootstrap approach. By this they mean an inductive argument for the $n$-bracket by \emph{solving} the equations for  the
 L$_\infty$ relations. They succeed for non-commutative Chern-Simons and for  non-commutative Yang-Mills 
 using properties of the \emph{star product} and various string field theories 
 combining Chern-Simons and Yang-Mills by careful computation
(see Section \ref{morphs} for uniqueness of their approach).
\p
 Contrast the inductive argument for BFV and BV (see Section \ref{bfv}) where a solution follows from the auxiliary acyclic resolution.

 \subsection{\LL-morphisms}\label{morphs}
 \h
 \LL-morphisms $V \to W$ can most efficiently be described as  morphisms of dgcas $C(V)\to C(W)$. They are particularly important in Kontsevich's proof of the
 Formality Conjecture in relation to deformation quantization of Poisson manifolds \cite{kont:defquant-pub}. Here the crucial \LL-morphism was from the dg Lie 
 algebra of multivector fields to the dg Lie of multidifferential operators. Since his initial success,
  there have been many formality theorems upgraded to  other settings, for example, for  BV algebras \cite{campos:BVformality} or  for a cyclic version \cite{ward:cyclicMC}.
   \p
 A special case is comprised of the Seiberg-Witten maps (SW maps) \cite{seiberg-witten} between ``non-commutative'' gauge field theories,  
 compatible with their gauge structures; in particular, between
 non-commutative versions of theories with a gauge freedom.  When these theories exist, they are consistent deformations of their commutative counterparts. 
 
\p
 Blumenhagen,  Brinkmann,  Kupriyanov and Traube \cite{bbkt}
establish uniqueness of their results (up to gauge equivalence) by constructing SW maps.
Aschieri and Deser \cite{deser-aschieri:SW} construct specific examples in terms of 
$U(n)$-vector bundles on two-dimensional tori given by
globally defined Seiberg-Witten maps (induced from the plane to the
torus). 

\subsection{BBvD \LL-algebras}\label{BBvD}
\h
In contrast to Lie $n$-algebras as defined above, it is possible to have  an \LL-algebra concentrated in degrees $0$ to $n-1$ with $d$ of degree $+1$
 as for the higher spin algebras of \BBvDs:
%The higher spin algebras of BBvD arise in  the following context:
They start with a given space of `fields' $\PP$ which is a module over  a Lie algebra 
$\XX$ of gauge symmetries.

%, which act as \emph{gauge transformations}. Although physically motivated, most of what follows does not depend on the nature of the  `fields',  the `gauge symmetries' or the operator.
\p

% on a given space of fields $\PP$.
 By a field dependent gauge transformation
 %\footnote{Mathematically speaking, we would say `an  action of $\XX$ on $\PP$', but we reserve `action' for the integrated Lagrangian.} 
 of $\XX$ on $\PP,$ they mean a polynomial (or power series)
map   $\Xi\otimes\LP\to \Phi$:
%$\Xi\otimes\Xi\oplus\Phi\to \Xi$:
\be
\dd _\xi(\pp) = T(\xi,\pp)=\Sigma_{i\geq 0} T_i(\xi,\pp)
\ee
 where
$T_i$ is linear in $\xi$ and polynomial of homogeneous degree $i$ in $\pp.$
Note the operation $T_0:\Xi \to \PP$ from `algebra' to `module', in contrast to  Lie $n$-algebras.
\p

They have a corresponding field dependent generalization of a Lie algebra structure on $\XX$:
a polynomial (or power series)
map $\Xi\otimes\Xi\otimes\LP\to \Xi$
% $\Xi\otimes\Phi\to \Phi$:
  $$[\xi,\eta](\pp) =  C(\xi,\eta,\pp)=\Sigma_{i\geq 0} C_i(\xi,\eta,\pp)$$
  where $C_i$  is bilinear in $\xi$ and $\eta$ and of homogeneous degree $i$ in $\pp.$
\p
 These operations obey consistency relations which Fulp, Lada and I  identified as  structure relations of an $L_\infty$ algebra \cite{fls}.

\begin{remark} Note the BBvD structure gives rise to an $L_\infty$-algebra structure on the direct sum of the space of fields and the space of gauge parameters, {\bf not} of the form of an \LL-algebra and its module, rather a Lie 2-algebra. This is similar to what occurs in the BFV and BV formalisms.
\end{remark}
\subsection{The BFV and BV dg algebra formalisms}\label{bfv}
\h For the related but separate notion of a BV algebra, see Definition \ref{bvalg}.
\p
Work of 
 Batalin-Fradkin-Fradkina-Vilkovisky (BFV) from 1975-1985 \cite{BF,FF,FV} concerned reduction of constrained Hamiltonian systems.
The constraints formed a Lie algebra and generated a commutative ideal. Their construction combined a Koszul-Tate resolution with a Chevalley-Eilenberg complex. 
The sum of the two differentials no longer squared to 0 but required terms of higher order. These were shown to exist using the acyclicity of the Koszul-Tate resolution.
 In 1988 \cite{jds:bull},   I was able to understand those terms in the context of homological perturbation theory (HPT, a common technique in $\infty$-theories)  and representations up to (strong) homotopy (RUTHs) (see Section \ref{ruths}).
\p
For the Lagrangian version, Batalin and Vilkovisky 
\cite{bv:closure, BV3,bv:anti} constructed the analogous dg algebra by a similar method \cite{jds:moskva} .
\p
These constructions  suggested an \LL-algebra and module, but things were not so simple.  
My student Lars Kjeseth in his UNC thesis \cite{lars:hha1,lars:hha2} showed that the appropriate structure was that of a \emph{strong homotopy} version of a  Lie-Rinehart algebra \cite{rinehart}. 
This concept then lay dormant until resurrected around 2013 in the work of Johannes Huebschmann \cite{jh:quasiLR} and of Luca Vitagliano \cite{luca:shLR}.
\p
More generally, \BFLS \ \cite{bfls:gauge} constructed an \LL-algebra on 
%A general construction of an sh Lie algebra ($L_{\infty}$-algebra) from 
any homological resolution of a Lie algebra. 
%The same construction applies for graded brackets infield theory such as the Batalin-Fradkin-Vilkovisky bracket of theHamiltonian BRST theory or the Batalin-Vilkovisky antibracket of the Lagrangian version.

%strong homotopy Lie algebras or L1-algebras.
\subsection{Some special  \LL-algebras}\label{dft}
\newcommand{\br}{[\cdot,\cdot ]}
\newcommand{\half}{\frac{1}{2}}

\h
Special \LL-algebras arise from Courant algebroids, Double Field Theory (DFT) and multisymplectic manifolds.
In the last two decades, what is called T-duality in string theory and supergravity led to a formulation of differential geometry on a generalized tangent bundle:
% (locally $TM \oplus T^*M$) as a \emph{Courant algebroid}. 

\begin{equation}
\label{eq:Edef}
   0 \longrightarrow T^*M \longrightarrow E 
      \stackrel{\pi}{\longrightarrow} TM \longrightarrow 0 . 
\end{equation}
Locally, the bundle $E$ looks like $TM\oplus T^*M$.
 \p
 \hskip10ex$\mathbf{Courant\  algebroids}$
 \h
Courant algebroids are structures which include as examples the doubles of Lie bialgebras and  bundles $TM\oplus T^*M.$ They are named for T. Courant \cite{courant}
who introduced them in his study of Dirac structures.
\p
%\begin{dfn}A \emph{Courant algebroid} is a sequnce of bundles $E\to TM \to M$ 
 %In fact, to  Lie 4-algebras where the graded pieces are.
%equipped with an operator $\DD: \cinf(M) \to \Gamma(E)$ as well as a bracket
%\end{dfn}

%\begin{dfn}
%$\label{def:jac} 
Given a bilinear skew-symmetric operation $[\ ,\  ]$ on a
vector space \( V \), its \emph{Jacobiator} \( J \) is the trilinear operator
on \( V \): 
\[
J(e_{1},e_{2},e_{3})=[[e_{1},e_{2}],e_{3}]+[[e_{2},e_{3}],e_{1}]+[[e_{3},e_{1}],e_{2}],\]
 \( e_{1},e_{2},e_{3}\in V \). 
%\end{dfn}
The Jacobiator is obviously skew-symmetric. Of course, in a Lie algebra \( J\equiv 0 \).

\begin{dfn}
\label{def:Liealgebroid} A \emph{Lie algebroid} is a sequence  $E\overset{\rho}{\to} TM \overset{p}{\to} M$ of vector bundle maps
with a Lie bracket on the space of sections $\Gamma(E)$
such that
\be
[X,fY] = \rho(X)(f)Y +f[X,Y]
\ee
for $X,Y\in\Gamma(E), f\in C^\infty(M). $ (It follows that $\rho[X,Y]=[\rho(X),\rho(Y)].$)
\end{dfn}

\begin{rmk}\label{rinehart} The purely algebraic analog of a Lie algebroid is a \emph{Lie-Rinehart algebra} over a commutative algebra more general than $C^\infty(M)$ \cite{rinehart}.
\end{rmk}

\label{def:quasi-algebroid} (Approximate Definition \cite{roytenberg:thesis,wein-royt:lmp})  A \emph{Courant algebroid} is an algebroid with anchor \( \rho :E\to TM \) 
%(in particular, a sequence of vector bundles $E\to TM \to M$) 
equipped with a nondegenerate symmetric  bilinear form  $\langle {\cdot }\,\ {\cdot } \rangle$
on the bundle $E,$ a skew-symmetric \emph{Courant}  bracket  $[\ ,\ ]$ on \( \Gamma (E) \) and a map \( {\mathcal{D}}:\cinf (M)\lon \Gamma (E) \)
satisfying many properties of which the most relevant is that the Courant brackets on  $\Gamma(E)$ satisfy the Jacobi identity up to a $\DD$-exact term:
 For any $e_{1},e_{2},e_{3}\in \Gamma (E) $,
\be
  J(e_{1},e_{2},e_{3})={\mathcal{D}}T(e_{1},e_{2},e_{3})
  \ee
  where \( T(e_{1},e_{2},e_{3}) \) is the function on the base \( M \) defined
by: 
\begin{equation}
\label{eq:T0}
T(e_{1},e_{2},e_{3})= 1/6  \underset{cyclic}{\sum}\langle [e_{1},e_{2}],e_{3}\rangle.
\end{equation}
\begin{remark}
A related structure is due to Dorfman \cite{dorfman:dirac} for which the bracket is not skew-symmetric but satisfies the 
 Loday (also called Leibniz)   version of the Jacobi identity \cite{loday} which is expressed as a derivation from the left
 (see Definition \ref{loday}). The name \emph{Dorfman bracket} appears in \cite{royt:dorfman} which discussed Courant-Dorfman structures.
 The Courant bracket is the skew symmetrization of the Dorfman  bracket \cite{yks:nij}.
\end{remark}
Roytenberg showed in his PhD thesis \cite{roytenberg:thesis} the equivalence of this 
structure with a specific $L_\infty$ algebra $X$, in which $l_1$ is determined by the de Rham differential, $l_2$ by the Lie bracket and the operation $l_3$ contains ``flux''-degrees of freedom (e.g. a three-form known as $H$-flux). Here $X$ is a resolution of ${\mathcal H}=coker\DD$:
\begin{equation}
\label{eqn:res}
X_2=\ker\DD\stackrel{d_2}{\lon}X_1=\cinf(M)\stackrel{d_1}{\lon}X_0=\Gamma(E)\lon{ \mathcal H}\lon 0,
\end{equation}
where
% $X_0=\Gamma(E)$, $X_1=\cinf(M)$, $X_2=\ker\DD$, 
with $d_1=\DD$ and $d_2$ is the 
inclusion $\iota:\ker\DD\hookrightarrow\cinf(M)$.
\p
 The Courant brackets on $ \mathcal H$ come from Courant brackets on  $\Gamma(E)$.
 %. for which the Jacobi identity is satisfied up to a $\DD$ exact term.
 %By they  appeal to  Theorem 7 in "The sh Lie Structure of Poisson Brackets in Field Theory" Barnich, Fulp, Lada, Stasheff) in 1998  constructs an \LL-algebra on a resolution of a space that has a Lie type bracket.   (cf. \cite{BFLS}), 
 Roytenberg and Weinstein \cite{wein-royt:lmp}
 use this to extend the Courant bracket to an
  \LL-structure on all of their   resolution $X$, manifestly a Lie 3-algebra.

%Following  Chris Rogers:
\p

 Further generalizations known as  \emph{higher Courant algebroids}  \cite{zambon:hicourant} are locally $TM\oplus \Lambda^k(T^*M)$.
\p
\hskip10ex $\mathbf{Double\  Field\  Theory (DFT)}$
\p
%As for Double Field Theory, 
Aldazabal, Marqu\'es and  N\'u\~nez write in
\emph{Double Field Theory: A Pedagogical Review}
\cite{ Aldazabal-Nunez}:

\begin{quote}Double Field Theory (DFT) is a proposal to incorporate  T-duality, a distinctive
symmetry of string theory, as a symmetry of a field  theory defined on a \emph{double configuration space}. 
\end{quote}

Indeed, as originally introduced in physics, DFT (Double Field Theory) refers \emph{not} to a doubling of fields but rather to doubling of an underlying structure such as double vector bundles and Drinfel'd doubles. The original Drinfel'd double occurred in the contexts of quantum groups and of Lie bialgebras. Following Roytenberg \cite{roytenberg:thesis},
Deser and I \cite{deser-jds} interpret the gauge algebra of DFT in terms of Poisson brackets on a suitable generalized Drinfel'd double. 
In DFT, it is the coordinates that are doubled and `double fields' refer to fields which depend on both sets of coordinates.
Thus we also refer to double functions, double vector fields, double forms, etc. In the physics literature, reference is made to a \emph{C-bracket} 
to distinguish it from a Poisson bracket,
but we showed that the C-bracket is essentially  the Poisson bracket on $T^*[1]TM$.
 (Here $[1]$ denotes the shift in degree of the \emph{fibre} coordinates of $TM$.)

\p
Similarly, Deser and S\" amann in \cite{deser-saemann:ERGI} used the analog of a Poisson bracket on $T^*[n]T[1]M$ to define an action of a Lie $n$-algebra
 on the algebra of functions 
$\mathcal{F}:= \mathcal{C}^\infty(T^*[n]T[1]M)$. In particular, what is called ``section condition'' in the physics of DFT 
was identified as the requirement for a specific Lie-2 algebra to act on $\mathcal F$ in a well-defined way.
\p
As a bundle over $TM$ with a local trivialization with respect to a covering $\mathcal U =\{U_\aa\}$, we have transition functions 
$a_{\alpha\beta}\in GL(2d,\mathbf R)$ satisfying the usual cocycle condition, but, in terms of the local splitting $TM\oplus T^*M$, there is a higher order `twist' $\tau$ depending on a 2-form $\omega$.  Now the cocycle condition fails, the failure depending on $\omega$:
  \[
g\dab\ g\dbgam \neq g\dagam \quad\mbox{on}\quad U\da\cap U\db\cap U\dgam .
  \]
This is often described as a failure of associativity, but it is more accurately
failure to correspond to  a representation. Since $g\dab,\  g\dbgam$ and  $g\dagam$ can be expressed in terms of 1-forms,  it can be that the difference is an exact form $d\lambda_{\aa\bb\gg}$
for some function $\lambda_{\aa\bb\gg}$ on $U\da\cap U\db\cap U\dgam.$ One would then say 
the transition functions form a representation up to homotopy (RUTH) or one can speak in terms of `gerbes'. (See Section \ref{wirth} for higher order generalizations.)
\p
\hskip10ex $\mathbf{Multisymplectic\  manifolds}$
\p
\begin{dfn} \cite{rogers:multisymplectic}
A  \emph{multisymplectic} or, more specifically \emph{$n$-plectic},  manifold is one
equipped with a closed nondegenerate differential form of
  degree $n+1$. 
  \end{dfn}
  As shown by Chris Rogers,
  \cite{rogers:multisymplectic}, 
  just as a
  symplectic manifold gives rise to a Poisson algebra of functions, any
  $n$-plectic manifold gives rise to a Lie $n$-algebra of differential forms with
  multi-brackets  specified via the $n$-plectic structure. The underlying graded vector space consists of a subspace of 
   $(n-1)$-forms he calls \emph{Hamiltonian} together with all  $p$-forms for $0 \leq p \leq
n-2:$
 \begin{equation} \label{complex}
\cinf(M) \stackrel{d}{\to} \Omega^{1}(M) \stackrel{d}{\to} \cdots \stackrel{d}{\to} 
\Omega^{n-2}(M) \stackrel{d}{\to} \Omega^{n-1}_{Ham}.
\end{equation}
%\subsection{Integration of \LL-algebroids}

% The bilinear bracket, as well as all higher $k$-ary brackets, are explicitly specified by the $n$-plectic structure.  

 \section{Generalized Jacobi identities and Nambu-Poisson algebras}
 \subsection{Identities}\label{n-ary}
 \h
There are two important ways to generalize to $n$-variables the Jacobi identity for a Lie algebra written as a  left derivation (see Definition \ref{loday}):
\be
[a,[b,c]]=   [[a,b],c] + [b,[a,c]].
\ee
In older literature, these generalized Jacobi relations are referred to as \emph{fundamental identities}, but, as suggested by de Azcarraga and Izquierdo \cite{Az&Iz:n-ary},
a better name might be \emph{characteristic identities}.,
\p
One identity  is the corresponding \LL-relation for a bracket of just $n$ variables:
\be
\sum \pm l_n(l_n(v_{\sigma (1)}\otimes \dots \otimes v_{\sigma (n}) \otimes v_{\sigma(j+1)}\otimes \dots\otimes v_{\sigma (2n-1)})=0. 
\ee

 It has   been studied quite independently of my work
and of each other by Hanlon and Wachs \cite{hanlon-wachs} (combinatorial
algebraists), by Gnedbaye \cite{gned} (of Loday's school) and by
de Azcarraga and Bueno \cite{az-bu} (physicists). 
\p
On the other hand,   
the characteristic identity for a \emph{Filippov algebra}  \cite{filippov})  says
$ [X_1, X_2, \dots,X_{n-1}, \quad]$ acts as a left derivation.

\begin{equation}
\label{n-der}
%\begin{aligned}
     [X_1, X_2, \dots,X_{n-1}, [ Y_1, Y_2, \dots, Y_n]]    =  
     \ee
     \be
   [[X_1,X_2,\dots, X_{n-1}, Y_1],Y_2,\dots,Y_n]     +  
 \dots + [Y_1,\dots, Y_{n-1},[X_1,\dots,X_{n-1},Y_n]].
%\end{aligned}
\end{equation}
That identity was known also to Sahoo and Valsakumar
\cite{sahoo}.
Unfortunately, both versions are called $n$-Lie algebras, hence my attempt (probably futile) to rename his as Filippov's.
\p
As I recall, I first learned of this other (Filippov) identity from Alexander Vinogradov when we met at the Conference on Secondary Calculus and Cohomological Physics, Moscow, August 1997. (See A. and  M. Vinogradov's
\cite{mmv} 
for a comparison of these two distinct generalizations of
the ordinary Jacobi identity to $n$-ary brackets.) The article  \cite{Az&Iz:n-ary} by de Azcarraga and Izquierdo
 is a very thorough survey of   even more $n$-ary algebras.
 \p
All these algebras  are important in geometry  and in physics
where the corresponding structures are on vector bundles over a smooth manifold (see \cite{wade} and references there in). 
 \subsection{Nambu-Poisson $n$-Hamiltonian mechanics}
 \p
Nambu's original work  \cite{nambu} was a generalization to an $n$-ary bracket of Hamiltonian mechanics with its binary Poisson bracket. Often the literature refers to Nambu-Poisson structures,
which emphasizes the setting of $C^\infty$ functions on a smooth manifold as in traditional  Hamiltonian mechanics. Just as the latter can be extended to sections of a Lie algebroid,
Nambu-Poisson structures on manifolds
can be  extended to the context of Lie algebroids
% in a natural way  based on the Vinogradov bracket 
 \cite{wade}. The Nambu bracket satisfies the Filippov identity \cite{takh:nambu}.
\p
One expects there to be ``$\infty$''-versions with the full panoply of applications as for \LL-structures, e.g. homological reduction of constrained  Nambu-Poisson algebras (cf. Section \ref{bfv}).
\p

 \section{Derived bracket and  brace and BV algebras}with major assistance from Fusun Akman, Yvette Kosmann-Schwarzbach, and Thedia Voronov
 \subsection{Differential operators}
 \h
An increasingly popular approach to \LL-structure is the use of \emph{derived brackets} (see Definition \ref{derivedbracket}) introduced by Kosmann-Schwarzbach   in \cite{yks:PtoG},  inspired by unpublished notes by Jean-Louis Koszul \cite{koszul:notes}.
 She traces their origin to work of 
Buttin and of A. Vinogradov in the context of  unification of brackets in differential geometry,
though  physics (of integrable systems for instance) was not far behind,  especially in work of Irina Dorfman \cite{dorfman:dirac,dorfman:book}. 
See  \cite{yks:derived} for an excellent survey and history.
\p
 Koszul's point of view  was that of  an algebraic characterization of the order of a differential operator. 
Multilinear operators with arguments that are left multiplication operators $\ell_{a}$ rather than elements $a$ of an algebra
 may  have appeared  first in various definitions of  the order  of   
 a   differential operator. For example, the \emph{order} of   a   differential operator $\D$ can be defined as
   \begin{dfn}\cite{akman:BV}
 \quad A linear operator $\D$ on an algebra $A$ is a \emph{differential operator of 
order $\leq r$} if an inductively defined $(r+1)$-linear form 
$\Phi_{\D}^{r+1}$ with values in $A$ is identically zero. 
  \end{dfn}
  or, more simply:
  \p
 \emph{An operator $\D$ on an algebra $A$ is of order at most $n$ if, for all left multipliers $\ell_a$ for $\in A$, the commutator $[\D,\ell_{a}]$
 % (this really needs to be the left multiplication operator, per Grothendieck's definition as well as to have both arguments on the same footing)}
  is of order at most $n-1$. }
\p

 This approach can be traced back to 1967:  Grothendieck (Ch.\,IV, EGA4), then developed by  Koszul~\cite{koszul:notes}.
 This was carried forward in 1997 by Akman \cite{akman:BV} who defines  higher order differential operators on a general 
\emph{noncommutative, nonassociative} graded algebra $\A$. Akman and  Ionescu \cite{akman:hiderived}  compare and show equivalence of several definitions of
order when the underlying algebra $A$ is classical, i.e.,  graded commutative and associative.
\p
The brackets that are now called  ``higher derived'' (as defined below) can be represented in a form similar to that for defining
``order'', but are of interest for  providing    \LL-structure. 
 \p
 \subsection{\LOO-algebras}
 \h
   Not only \LL-algebras are relevant, but also \emph{\LO}-algebras and even \LOO-algebras. Loday (1946 - 2012) originally called them \emph{Leibniz algebras}, but he  deserves the credit.
  Recall:
  \begin{dfn}\label{loday}
  A (left) \emph{\LO-algebra} $(V,[\ ,\ ])$ consists of a vector space $V$ with a \emph{Loday bracket} $[\ ,\ ]$ satisfying the left Leibnitz identity:
  \be
  [a,[b,c]]= [[a,b],c] \pm [b,[a,c]].
  \ee
    \end{dfn} 

\begin{remark} Throughout this section, ``there exists a set of signs'; here for the graded case.
\end{remark} 

 \begin{dfn}\cite{loday:coprod}
  A (left) \emph{\LOO-algebra} $(V,\pi)$ consists of a graded vector space $V$ with a 
  sequence $\pi=(\pi_1,\pi_2,\dots)$ of multivariable coderivations $\pi_i$ on the free graded coalgebra cogenerated by 
  $V$ shifted with a rather unusual coproduct using unshuffles that preserve the (last) element. 
  The sequence $\pi$ is to satisfy  $[\pi,\pi]=0$.
  \end{dfn} 

Apparently this was defined first by Livernet   \cite{livernet:lodayinfty} for right \LOO-algebras, later independently in a different setting by Ammar and Poncin \cite{ammar-poncin:lodayinfty}
with further work by Peddie \cite{peddie:lodayinfty}.

\subsection{Derived brackets}\label{derived}
Kosmann-Schwarzbach  defined \emph{derived brackets} in 1996 \cite{yks:PtoG}:

 \begin{dfn}\label{derivedbracket}\cite{yks:PtoG} If $(V, [~,~], D)$ is a graded 
differential Lie or Loday algebra
%(Yvette's original construction) algebra over $R$ 
with a bracket $[\ ,\ ]$of degree $n$,
the {\em derived bracket} of $[~,~]$ by $D$ is  denoted 
$$[~,~]_D: V \otimes V \to V$$ 
and defined  by
\be
[a,b]_D= (-1)^{|a|+n+1} [Da,b] \ ,
\ee
for $a$ and $b \in V$.
\end{dfn}

Often $D$ is itself  given as $D=[\D,\  \ ]$ for an element $\D\in V$.
\p
In this generality, the derived bracket is graded Loday and respected by the differential. To obtain a differential graded Lie bracket, 
additional actions are necessary.
\p

Continuing the iteration of Lie brackets led  in 1996 to Bering  defining  \emph{higher BV anti-brackets}  in physics  \cite{beringetal:hibrack}.
Somewhat later, independently in math, T. Voronov defined   \emph{higher derived brackets} \cite{tv:hiderived}.  This was followed quickly by \cite{bering:brackets} 
and further proliferation. Apparently the membrane between math and physics is osmotic: information travels one way faster than the other. 
\p
The definition closest to my interests is:
\begin{dfn} \emph{Higher derived brackets}
\[ B^r_{\Delta}(a_1,\dots,a_r)\] are defined as
\[ [\dots [\,[\Delta,a_1],a_2],\dots,a_r]\]
where $(L,[-,-])$ is a Lie algebra with $\D, a_1,\dots,a_r\in L$. 
\end{dfn}
In terms of operators $D$:
\begin{dfn} [Alternative] \emph{Higher derived brackets}
\[ B^r_{D}(a_1,\dots,a_r)\] are defined as
\[ [\dots [Da_1],a_2],\dots,a_r]\]
where $(L,[-,-])$ is a Lie algebra with $ a_1,\dots,a_r\in L$ and $D:L\to L$. 
\end{dfn}
The major uses of higher derived brackets are for describing/constructing  \LL-algebras and other $\infty$-algebras. An important case is the derived bracket 
built with an odd, square-zero  differential operator   $\D$ of order $\leq 2$, which Akman uses
to define a \emph{generalized Batalin-Vilkovisky algebra} structure on an algebra $A$ (see  Section \ref{bvalg}). The emphasis on  order $\leq 2$ is historical,
 dating back to \BV\  and their applications in physics..
\p
%To be added in a later version: Discussion of abelian subalgebra case \cite{tv:hiderived} and more general \cite{bandiera:hiderived}.
\p
T. Voronov's   \emph{higher derived brackets} \cite{tv:hiderived} first occurred in the context of a  Lie algebra $L$ with an abelian sub-algebra $A$.  The 
higher derived brackets on $A$ were derived on $L$ and then projected back to $A$. More general higher derived brackets are introduced by Bandiera  
\cite{bandiera:hiderived}.
% \cite{akman:BV}.
% as well as the $\infty$ version (see Section \ref{gbv}).
\p
Derived (binary and higher) brackets are also important for describing and understanding infinitesimal symmetry actions relevant in physics.
The Roytenberg-Weinstein \LL-structure can be expressed and generalized in terms of  derived brackets \cite{roytenberg:thesis},
Most recently, Deser and S\"amann \cite{deser-saemann} show binary derived brackets underlie the symmetries of Double Field Theory (Section \ref{dft}).
They suggest adopting as a guiding principle: 
\begin{quote}Whenever we are seeking an infinitesimal action of a fundamental symmetry, we try to find the corresponding derived bracket description. 
\end{quote}
One should also  seek higher derived brackets as in \cite{tv:hiderived} and
hence  \LL-structures as called for in \cite{hz} (see Section \ref{LLphysics}). 

\subsection{Braces}
Brace algebras, like operads, are all about composition. As Akman says: it all depends on what you want to achieve and how far you want to go.
\begin{dfn} \cite{gv:hiops, B&G}
A \emph{brace algebra} is a graded vector space with a collection of braces
$x\{x_1, \dots, x_n\}$ of degree $\pm n$ satisfying the identities 
\begin{multline}
\label{higher}
x\{x_1, \dots, x_m\} \{ y_{1}, \dots , y_{n}\} \\
\begin{split}
= \sum
\pm x\{ y_1, \cdots,
y_{i_1}, x_1  \{ y_{i_1+1}, \cdots,
 y_{j_1}\}, y_{j_1+1}, \cdots \hskip10ex \\
 \cdots,
y_{i_{m}}, x_m \{  y_{i_{m}+1} ,
\dots , y_{j_m}\}, y_{j_m+1}, \dots, y_n \} 
\end{split}
\end{multline}
where the sum is over $0 \le i_1 \le \dots \le i_m \le n$ and
the sign is due to the $x_i$'s passing through the $y_j$'s in the
shuffle.
\end{dfn}

Thus brace algebras are flexible enough to accommodate insertion of several arguments at one time, but not necessarily filling all slots (an operad does fill all slots).
\p
Braces first occurred in work of Kadeishvili \cite{kad:braces} who called them higher $\smile_1$ products. In later but independent work,  
Gerstenhaber and A. Voronov \cite{gv:hiops, B&G} named  them \emph{braces}, which is now the common term used by several authors.
 There is a corresponding \emph{brace operad} with subtle use of  trees; \cite{dolgushev-willwacher:braceoperad} explains its relation to Deligne's `conjecture' for \AA-algebras.
 \p
(Higher) derived brackets  can be defined in terms of braces (special compositions of maps).
\p

\subsection{BV algebras}\label{bvalg}

%\vskip2ex

\dbend \ {\bf Warning!!} 
\p
Batalin and Vilkovisky have made two \emph{distinct} contributions to gauge field theory:
\begin{itemize}
\item the formalism of the Lagrangian dg construction with  antifields and ghosts etc. (see \ref{bfv}).
\item graded algebras with a special  differential operator of order 2.
\end{itemize}
Unfortunately both are sometimes refered to as the BV \emph{formalism}, which preferably should be used for their dg construction.
\p

Recall a \emph{Gerstenhaber algebra} $(A,\cdot,[\ ,\ ])$ is the `odd' analog of a graded Poisson algebra, $[\ , \ ]$ being called a Gerstenhaber bracket. 
The first example 
occured in his paper on Hochschild cohomology of an associative algebra \cite{gerst:coh}.
\begin{dfn}
A \emph{Batalin-Vilkovisky algebra}, or \emph{BV algebra}, is a Gerstenhaber
algebra where the bracket $[ \; , \; ]$ is obtained from an odd, 
square zero, second order differential operator $\D$:
\be
 [a,b]  =(-1)^{|a|}\D(ab)-(-1)^{|a|}\D(a)b-a\D(b) 
 \ee
\end{dfn} 
Batalin and Vilkovisky introduced the algebras now named in their honor in their study of conformal field theroies and quantization \cite{bv:quant}.
\p
There is a notion of \emph{generalized BV algebra} due to Akman \cite{akman:BV} in which the bracket $[ \; , \; ]$ is obtained from an odd, 
square zero, second order differential operator $\D$ as above, though  the algebra doesn't have to be associative nor commutative,
but   the  bracket still measures the deviation of $\D$ from being a derivation. She also considers 
\emph{differential BV algebras}. For Akman, a major motivation is the work of Lian-Zuckerman \cite{lz:new} related to VOAs (vertex operator algebras).
The restriction to `second order' is just the classical \BV\ case.
\p

Of course, there is a `higher' notion of $BV_\infty$ algebra. The article by G\'alvez-Carrillo,   Tonks and Vallette \cite{vallette-etal:homotopyBV}  should be the canonical reference.
They point out  that several authors call a ``homotopy BV-algebra'' what is a \emph{commutative} homotopy BV-algebra, that is one for which the product is commutative
(compare remarks at the end of Section \ref{iterated}).
In addition to  giving an  operadic description (see also \cite{drummond-cole-vallette}), they give four equivalent definitions of a homotopy BV-algebra and provide  applications in four different categories. 
In particular, they develop the deformation theory and homotopy theory of BV-algebras and of $BV_\infty$ algebras. 

\section{Deligne's question and cyclicity}\label{quest}
\h In his seminal work \cite{gerst:coh}, Gerstenhaber showed the Hochschild cohomology $HH^\bullet(A)$ of an associative algebra 
carries what is now known as a Gerstenhaber algebra structure.   Moreover, he constructed homotopies on the Hochschild cochains $CH^\bullet(A)$
to yield the relevant identities on $HH^\bullet(A)$. Thirty years later
 in 1993 in a letter to several of us, Deligne \emph{asked}: Does the Hochschild cochain complex
of an associative ring have a natural action by  chains of the small
squares operad?
\p 
%The first few low level cochain homotopies were present in Gerstenhaber.
 This lead to \emph{many} `higher structure' papers resolving this `conjecture'  and generalizations,
  in many of which  the little disks operad was invoked instead. In some, braces play a key role.
  The chains involved were various singular chains on the topological operad or PROP or  cellular chains on a related CW complex.
  The generalization to the Hochschild cochains when $A$ is an \AA-algebra followed soon after.
  To follow this development further means entering a labyrinth and trying to choose a path (see \cite{mcclure-smith,kaufmann:spineless} among others). That of \cite{MandY:defoveroperad} emphasizes relations to deformation theory (see Section \ref{defs}).
 \p
  Let me instead pay attention to  the cyclic Deligne `conjecture'  which, on the one hand,   refers to the Hochschild cochain complex for an algebra $A$ with an invariant inner product and, on the other hand, to the \emph{framed little disks operad} \cite{getz:bv2d}. Solutions   involve \emph{cacti} \cite{sasha:cacti,kaufmann:cacti},  \emph{spineless cacti} \cite{kaufmann:spineless} and Sullivan chord diagrams \cite{tz:aso,tz:cyclic} and of course compactified moduli spaces. The solution is developed further in an \AA\  context and   treated particularly well  by Ben Ward  \cite{ward:cyclic}.
Although the latter clearly involves associativity and cyclicity, it is intriguing to see  the approach of Kaufmann and Schwell \cite{kaufmann-schwell} involves both the associahedron and the cyclohedron, which was introduces by Bott and Taubes \cite{BT} though only later named in \cite{jds:config}. Progressing one step further, Tradler considers an $\infty$-inner product \cite{tradler:inftyinner}.
   
\section{\LL\  in deformation theory}\label{defs}
\h
Deformations of complex structure go back to Riemann, but just for one complex dimension. In higher dimension, even  a proper definition of infinitesimal deformation
was a problem, resolved by Nijenhuis and Fr\" olicher \cite{fn} by identifying it as an class in the cohomology $H^1$ of the sheaf of germs of holomorphic tangent vectors. 
Kodaira and Spencer took over from there \cite{ks:I&II}. The primary obstruction to extending the infinitesimal $\delta$ was a class in $H^2$ given by the product 
$\delta\smile\delta$.
Higher obstructions drew the attention of Douady \cite{douady}, who related them to Massey `products' . Unnoticed for a while, this approach was carried further by Retakh \cite{retakh:massey}).
\p
\subsection{What is deformation theory?}
\h  Based on the history  for complex structures, Gerstenhaber offered the first general description:
\begin{quote}
a deformation theory seems to have at least the following aspects:
\begin{itemize}
\item A definition of the class of objects within which deformation takes place, and identification of the infinitesimal deformations of a given object with the elements of a suitable cohomology group.
\item A theory of the obstructions to the integration of an infinitesimal deformation.
\item A parameterization of the set of objects obtainable by deformation from a fixed one, and the construction of a fiber space over this space, the fibers of which are the objects.
\item A determination of the natural automorphisms of the parameter space (the modular group of the theory) and determination of the rigid objects. In some cases almost all points of a
 parameter space will represent the same rigid object, degenerating in various ways to objects admitting proper deformations.

\end{itemize}
\end{quote}
%Closely related to Gerstenhaber's description is one in terms  of associated \emph{moduli spaces}, as in Riemann's original work.
By analogy with Riemann's original work, the  associated parameter spaces are referred to as \emph{moduli spaces}. In contrast to classical results, the moduli space need not be in the class of objects being studied.
\p
%The algebraic version dates back to the work of Murray Gerstenhabe% (1963 - it was a very good year!!),who 
 From  1963 to 1968,  Gerstenhaber studied deformations of associative algebras in terms of Hochschild cochains and cohomology \cite{gerst:defs-ra}. 
 This had  major impacts: for deformation theory, for physics (see Section \ref{quant}) and for higher homotopy theory (see Section \ref{quest}). His initial work was soon followed by work of Nijenhuis and Richardson 
 ,\cite{nij-rich:algdef, nij-rich:Liedef}. 
They  were the first to articulate something close to the current ``metatheorem'', adding dg Lie algebra as a crucial ingredient: 
\p
\emph{The deformation theory of any mathematical object, e.g., an associative algebra, a complex manifold, etc., can be described starting from a certain differential graded (dg) Lie algebra associated to the mathematical object in question.}\label{meta} 
\p
This provides  a natural setting in which to pursue
the obstruction method
for trying to integrate ``infinitesimal deformations''. The deformation equation is
known as the Master Equation in the physics and physics inspired literature and now most
commonly as the Maurer-Cartan equation. It is also the equation for a \emph{twisting cochain} used in describing twisted tensor products \cite{gugjds,jds:twisted}.
\p
As in the above examples, some authors refer to the corresponding \emph{cohomology} as controlling the deformations.
\p
In 1986, inspired instead by Goldman-Millson \cite{goldman-millson}, Deligne stated this as a philosophy in a letter to Millson \cite{deligne:gm}.
%: https://publications.ias.edu/deligne/paper/2595.
\p
Meanwhile in 1967, Lichtenbsum and Schlessinger \cite{schlessinger-lichtenbaum} expressed deformation theory in terms of a \emph{cotangent complex}; think differential forms on the moduli space with 
the \emph{tangent complex} corresponding to a Lie algebra of derivations. This  was carried further in Mike's thesis \cite{schlessinger}, written from the viewpoint of representability of certain functors.
\p
The extension of deformation theory to a variety of mathematical objects is well surveyed by Gerstenhaber and Schack \cite{gerstenhaber-schack:survey}.

%\trd{??mention paper that uses simplicial Artin rings??}
 \subsection{The (not quite) metatheorem}\label{quite}
 \p
 In his MR  review of \cite{MandY:defoveroperad}. A. Voronov writes with respect to the `metatheorem':
 \begin{quote}
 The uncomfortable generality of the statement might have prevented anybody from making it a theorem and proving it. A route chosen by many was usually the following one: given a mathematical object $A$, construct a dg Lie algebra $L$ and answer as many questions about the deformation theory of $A$ in terms of $L$ as you can.
 \end{quote}
 \p
 %Closely related to Gerstenhaber's description is one in terms  of associated \emph{moduli spaces}, as in Riemann's original work.
  There have been  several theorems of great but not complete generality.
 \p
  Markl \cite{markl:cotan} develops a deformation theory phrased as a controlling cohomology theory for $k$-algebras over  $k$-linear equationally given category, e.g. over operads and PROPs - for bialgebras.. 
He invites comparison to traditional cotangent cohomology of a commutative algebra (in characteristic zero) based on the free differential graded
algebra resolution of the algebra under consideration \cite{Sjds, schlessinger}. In
% $\emph{Intrinsic brackets and the ${L_\infty}$-deformation theory of bialgebras} 
\cite{markl:JHRS10}, he provides explicit constructions in terms of ${L_\infty}$-deformation theory, enhanced further in \cite{markletal:algdiags}; 
see an alternative  \cite{merkulov-valletteI} for
representations of prop(erad)s.
\p
Others are  expressed in `languages the muse did not sing at my cradle'.
Kontsevich and Soibleman \cite{MandY:defoveroperad},  develop the existence of deformation theory for an algebra over any (colored) operad, but in  
a more geometric language of formal dg manifolds. This refers to  cofree cocommutative coalgebras with differential in disguise. Thus again \LL-algebras are doing the controlling.
They give further insight in the comparison of the algebriac and geometric version of \AA-structures in \cite{kont-soi:notes}.The one closest to my native language e.g. of \LL-algebras, is that of Pridham \cite{pridham:unifying,pridham:corrig}, which is the most general, going beyond characteristic 0. 
Closer to (derived) algebraic geometry is that of Lurie \cite{lurie}.

\subsection{Deformation quantization}\label{quant}
\h
Especially important was the eponymous \emph{Gerstenhaber bracket}  \cite{gerst:coh} which was later identified by interpreting the Hochschild complex in terms of coderivations \cite{jds:intrinsic, akman:chicken}. His work led to an algebraic description of  \emph{deformation quantization} \cite{bfflsI,bfflsII}, a term   derived from physics.
%, but really refers to an algebraic issue:
\p
Given a Poisson algebra $(A, \{\ ,\  \})$, a deformation quantization is an associative unital $\star$ product  on the algebra of formal power series  $A[[\hbar]]$ subject to the following two axioms:
\begin{align}
 f\star g &=fg+\mathcal{O}(\hbar)\\
 f\star g-g\star f & =\hbar\{f,g\}+\mathcal{O}(\hbar^2).
\end{align}

In the physics context, the Poisson bracket is given in terms of differential operators as are the sought after  terms of higher order (in which the devil resides!).

\subsection{\LL-algebras  in rational homotopy  theory}\label{qht}
\h
\LL-algebras    arose by 1977 in my work with Mike Schlesinger \cite{Sjds,SS, jds-ms}  on deformation theory of rational homotopy types. %as well as in mathematical physics. . 
%The yoga of deformation theory, that any problem in deformation theory is ``controlled" by a differential graded Lie algebra (unique up to homology equivalence of dg Lie algebras, which is to say: via an \LL-morphism) (see Section \ref{morphs}).
 Mike and I extended the yoga of control by a dg Lie algebra  (See\ref{meta}) to similar control by an \LL-algebra,
for example, on the homology of a huge strict dg Lie algebra.
\p
\newcommand{\HH}{\mathcal H} 
 Let $\HH$ be a simply connected graded commutative algebra of
finite type and $(\Lambda Z,d) \to \HH$ a filtered model.
Differential graded Lie algebras provide a natural setting in which to pursue 
the obstruction method
for trying to integrate ``infinitesimal deformations'', elements of
$H^1({Der} \Lambda Z)$, to full perturbations.  In that regard,
$H^*({Der} \Lambda Z)$ appears not only as a graded Lie algebra (in
the obvious way) but also as an \LL-algebra. 
 \p   
Our main result compares  the \emph{moduli space} set of augmented
homotopy types of dgca's 
$(A,i : \HH \approx H(A))$ with the path components of $C(L)$ where $L$ is a sub Lie algebra of 
$ {Der}\,\Lambda Z$ consisting of the weight
decreasing derivations.
We were inspired by the  work of Kadeishvili on the \AA-structure on the homology of a dg associative algebra. This 
again used the `higher structure' machinery of HPT (Homological Perturbation Theory). 

\p
Independently, in 1988,  correspondence between  Drinfel'd and Schechtman  \cite{Drinfeld:letter, Toen:drinfeld} 
develops \LL-algebra under the name \emph{Sugawara - Lie algebra} for the needs of deformations theory.
%\footnote{
 As with `Jacobi up to homotopy', the need for \LL-algebras was `in the air' (See Section \ref{Loo}). 
 ``A letter from Kharkov to Moscow'' \cite{Drinfeld:letter} deserves further attention. 
 \p
 The disparity  in the dates recalls the lack of communication between the fSU and the West. The situation was much better two years later 
 when Gerstenhaber and I were able to participate in person at the Euler International Mathematical
Institute's
Workshop on {Q}uantum {G}roups,
Deformation Theory, and Representation Theory \cite{leningrad1990}. Interactions there were to prove very fruitful (see Section \ref{fields}).

%This led to his \cite{kont:trium}.
%- see Section \ref{cyclic}}.
\section{Representations up to homotopy/RUTHs}\label{ruths}
\h
A good bit of group theory has been carried over to \AA-algebras,  but far from  all.  Most recently, the analog of Sylow theorems appeared \cite{prasma:sylow}.
It was late in the game before higher homotopy representations received much attention. Although the language is slightly different, an ordinary representation of a group or algebra is equivalent to a morphism to the endomorphisms of another object. This appeared early on in homotopy theory (I learned it as a grad student from Hilton's  \emph{Introduction to Homotopy Theory} \cite{hilton:intro}  - the earliest textbook on the topic) in terms of the action of the based loop space $\Omega B$ on the fiber $F$ of a fibration $F\to E\to B$:
$$\Omega B \times F\to F \ \  or \ \  \Omega B \to Haut(F)$$
where $Haut(F)$ denotes the monoid of self homotopy equivalences $F\to F$.
The existence of $\Omega B \times F\to F$ follows from the covering homotopy property, but with uniqueness only `up to homotopy'.
\p
 Initially, this was referred to as a \emph{homotopy action},
meaning only that $(f \lambda)\mu$ was homotopic to $f(\lambda)f(\mu),$ with no higher  structure. 
The map $\Omega B \to Haut(F)$ is  only   an H-map; the full equivalence between such fibrations and \AA-maps $\Omega B \to  Haut(F)$ had to wait until such maps were available.
%; here $Haut(F)$ denotes the self homotopy equivalences $F\to F$. 
The corresponding terminology is that of (strong or $\infty$) homotopy action, which has further variants under a variety of names. The definition is simplified if we use the Moore space of loops  $\MB$ \cite{moore:loops} (see Section \ref{mooreloops}).

\begin{definition}
\label{Thetas}
A \emph{representation up to homotopy} of  $\MB$ on a fibration $E\to B$ is 
an $A_\infty$-morphism (or shm-morphism \cite{sugawara:hc})  from $\MB$ to $End_B(E)$ (see Definition \ref{shmap}).
\end{definition}
 
In 1971,
Nowlan  \cite{nowlan:action} considered fibrations $F\to E\to B$ with associative H-space F as fibre,
  acting fibrewise on E, but the action being  associative only up to higher homotopies. His main result  shows  that a fibration 
$\Omega Y \to E \to B$ is fibre homotopy equivalent to one induced by a map of $B \to  Y$ if and only if $E$ admits an \AA-action of  $\Omega Y$.
% (An alternative solution to the problem also involving higher homotopies has been given by G. J. Porter [Illinois J. Math. 16 (1972), 41Ð60].) 
\p
\AA-actions occur also for relative loop spaces \cite{hoefel:Aooaction,vieira:rel}.
\p
In the case of a smooth vector bundle $E \to B$, the corresponding notion is parallel transport, usually determined by a \emph{choice} of connection, hence unique.
In physics, various action functionals for quantum field theories correspond to  \emph{higher parallel transport} in bundles with graded vector space fibers, corresponding
to the higher homotopies of the strong homotopy action.
%From a topological point of view, strong homotopy action is the more basic notion, not depending on a connection but not unique except up to appropriate higher homotopies. 
%\p
%A Poisson algebra is a commutative associative algebra A with an (anticom- mutative) bracket { , } which is a derivation with respect to the commutative product:
%$\{f,gh\} = \{f,g\}h + f\{g,h\}$. Constraints constitute a distinguished set of elements (j)a of A. They are said to be first class constraints if the ideal ? they generate (under the commutative product) is closed under Poisson bracket; I need not be an ideal with respect to { , }. This structure arises in physics with A = C°°(W) for some symplectic manifold W
 
   \p
   Switching to algebra, there is the corresponding notion of \emph{representation up to homotopy} of an associative dg algebra on a dg vector space $V$.   This seems to have first occurred in \cite{jds:hrcpa} in the context of a Poisson algebra $(P,  \{ \ ,\ \})$ with a commutative ideal $I$ closed under $ \{ \ ,\ \}$. Such a structure arises in physics with 
   $A = C^\infty(W)$ for some symplectic manifold, $W$ (see Section \ref{bfv} \cite{BF,FF,FV}).
   \p
   There are analogous notions of RUTHs of Lie structures. In particular, in open-closed string field theory  (see Section \ref{ocha}), the CSFT acts up to homotopy in the strong sense on the OSFT and similalry for OCHAs with the \LL-algebra acting on the \AA-algebra by $\infty$-homotopy derivations..
   In 2011, RUTHs of Lie algebroids was developed by Abad and Crainic \cite{abad-crainic:shreps}\footnote{As of this writing, the wiki is not up to date.}. In particular, they use representations up to homotopy to define the adjoint representation of a Lie algebroid to  control  deformations of  structure (see Section \ref{defs}).
 
   \section{\AA-functors, \AA-categories and $\infty$-geometry}\label{wirth}
   \h
  Since associativity is a key property of categories, it is not surprising that  \AA-categories were eventually defined. In 1993, Fukaya \cite{fukaya} defined them to handle  Morse theoretic homology.
 \p 
  Just as one considers \AA-morphisms of \AA-algebras, one can consider \AA-functors 
  (also known as \emph{homotopy coherent functors}) between \AA-categories. Such functors were  first considered for ordinary strict but topological categories in the context of classification of fibre spaces. For fibrations which are   locally
homotopy trivial  with respect to a good  open cover $\{ U\da\}$ of the base, one can define transition functions
  \[
g\dab : U\da\cap U\db \to  Haut(F),
  \]
%where $H$ is the monoid of homotopy equivalences of $F$ to itself
 but instead of the cocycle
condition for fibre bundles, one obtains only that $g\dab g\dbgam$ is homotopic to $g\dagam$ as a map of
$U\da\cap
U\db\cap U\dgam$ into $Haut(F)$. 
\p
 In 1965, Wirth \cite{wirth:diss, jw-jds} showed how a set of coherent higher homotopies arise on multiple intersections. He calls that set a \emph{homotopy transition cocycle.}
The disjoint union
$\coprod U\da$ can be given a rather innocuous structure of a topological category $U$,
i.e., $\Ob U = \coprod U\da$ and $\Mor U=\coprod U\da\cap U\db$.
%that is $x\circ y=x=y$ is defined iff $x\in U\da$, $y\in U\db$ and $x=y$.  
Regarding 
$Haut(F)$ as a category with one object in the standard way, Wirth shows 
the transition cocycle  web of higher homotopies
is precisely equivalent to a  homotopy coherent functor.
For `good' spaces, the usual classification of such fibrations is effected by the realization of this functor via a map $BU\to BHaut(F).$
%\subsection{$\infty$-geometry?}
\p
Having come this far in the \AA world, what about \emph{`$\infty$-geometry'}?
\p
Recall Roytenberg and others define a \emph{dg manifold} as a locally ringed space
 $ (M, \mathcal{O}_M)$  (in dg commutative algebras over $\mathbb R$), which is locally isomorphic to $(U,  \mathcal{O}_U )$, where
$ \mathcal{O}_U  = C^\infty(U) \otimes S(V^\bullet)$ where $\{U\}$ is an open cover of $M$, $V^\bullet$ is a
 dg vector space and $S(V^\bullet)$ is the free graded commutative algebra.  
 Again the  transition functions $g\dab: U\da\cap U\db \to Aut S(V^\bullet)$ with respect to an open cover $\{ U\da\}$ of $M$ satisfy  the cocycle condition.
\p
 Notice that `manifold' is irrelevant, a hold-over from the early days, but useful intuition for those trained in differential geometry; at most a topological space with a good open cover is needed. The coherent homotopy generalization of the definition of a dg-manifold is straightforward, but requires a  coherent homotopy cocycle condition, as proposed here:
 \vskip2ex
 \begin{dfn} A \emph{dg $\infty$-manifold} or \emph{sh-manifold} is a locally ringed space
 $ (M, \mathcal{O}_M)$  (in dg commutative algebras over $\mathbb R$), which is \emph{locally} homotopy equivalent (as dcga's)  to ($U,  \mathcal{O}_U )$, where
$ \mathcal{O}_U  = C^\infty(U) \otimes S(V^\bullet)$ with $\{U\}$ an open cover of $M$ and 
$(S(V^\bullet),d)$ is free as graded commutative algebra. 
 \end{dfn}
 \p
The higher analogs of classical transition functions with a cocycle condition are exactly Wirth's   homotopy transition cocycles.

\section{\AA-structures  and  \LL-structures in physics}\label{aall}
\h
There has not been  as much presence  of \AA-structures   in physics as of as of \LL-structure, though Zwiebach followed with \AA-algebras for open string field theory in 1997 \cite{G&Z:open},
a few years after his \LL-algebra for closed string field theory.  Then,  for open-closed string field theory, \AA and \LL were combined.
\p
On the other hand, discovery by physicists of `mirror symmetry'  among Hodge numbers of dual Calabi-Yau
manifolds led Kontsevich to propose  \emph{homological mirror symmetry} \cite{kont:hms} in terms of \AA-categories such as were presented by Fukaya \cite{fukaya} (see Section \ref{wirth}).
His suggestions included extended moduli spaces  and deformations of \AA-categories.
%of the total Hochschild cohomology of an ?extended moduli space? $M$
\subsection{String  algebra}
\p
In OSFT, strings are considered as paths in a manifold with  interactions handled by "joining" two strings to form 
a third. This is often pictured in one of three ways (see \cite{jds:almost} for graphics): 
\begin{itemize}
\item E:  endpoint interaction:
(most familiar in mathematics as far back as the 
 study of the fundamental group) occurs only when the end of one string agrees with the beginning of the other and the parameterization is adjusted appropriately,
 \item M:  midpoint or half overlap interaction: occurs only when the last half of one agrees with the first half of the other as parameterized,
 \item  V:  variable overlap interaction: occurs only when, as parametrized,  a last portion  of one agrees with the  corresponding first portion of the other.
 
 \end{itemize}
\p
The endpoint case E comes in two flavors:
\begin{itemize}
\item with paths parameterized by a fixed interval $[0,1]$ in math or $[0,\pi]$ or $[0,2\pi]$ in physics,
\item with paths parameterized by intervals $[0,r]$ for $r\geq 0$ (see Section \ref{mooreloops}).
\end{itemize}
\p
The midpoint case  M  was considered 
by Lashof \cite{lashof:loops} and later by Witten \cite{witten:noncommsft}; it has the advantage of being associative when defined.  Note for the iterated composite of 
three strings the middle one disappears in the composite! For Lashof and topologists generally, the interaction was regarded as $a+b=c$ whereas for  Witten and physicists 
generally treat the interaction of three 
 strings more symmetrically, cf.  by reversing the orientation of the 
 "composite'': $a+b= -c$ or $a+b+c=0$. One could call this a `cyclic associative algebra'.   According to Kontsevich, his reading of 
  \cite{jds:almost} led to defining \emph{cyclic \AA-algebras}.
  % (see Section \ref{fields}).
\p
The variable case  V  seems to have occurred first in
physics in the work of Kaku \cite{kaku}.  Surprisingly, it is associative only up to homotopy but the pentagon relation holds on the nose! 
In physspeak, that is said as 3-string and 4-string vertices suffice \cite{hikko}.
\p
The images of closed string interaction and open- closed string interaction are much more subtle.

\subsection{String  field theories}\label{fields}
\h   
 But those are for strings; string fields are functions or forms on 
the space of strings and they form an algebra under the convolution product
where the comultiplication on a string is the set of \emph{decompositions } at  arbitrary points in the parameterizing interval. 
(An excellent and extensive `bilingual' (math and physics) survey is given by Kajiura \cite{kajiura:survey2010}.)
\p
One of the striking  aspects of \AA- and \LL-algebras in physics is the use of an inner product 
$<\ ,\ >$ for the \emph{action functional}, an integral of a real or complex valued function
of the fields.
I first noticed this in Zwiebach's CSFT \cite{z:csft} for the classical (genus 0) action 
\be
\sum\int (\phi_0, \phi_1,\cdots,\phi_n)
\ee
where the $\phi$'s are string fields and the integrand is cyclically symmetric (up to sign).
 \p
 The corresponding \LL-structure is  determined by
\be
<v_0,\ l_n(v_1,\cdots,v_n)>\  =\  (v_0, v_1,\cdots,v_n).
\ee
Cyclic \LL-algebras were formalized by Penkava \cite{penkava:cyclic}.

\p
As best I can determine, it was Kontsevich who first considered \AA-algebras with an invariant inner product, know as  \emph{cyclic \AA-algebras}.
Retakh recalls he and Feigin  discussed one of  the talks at the Euler Workshop on {Q}uantum {G}roups,
Deformation Theory, and Representation Theory and a text associated with the talk and showed it to Kontsevich or as Kontsevich says:
\begin{quote}
B. Feigin, V. Retakh and I had tried to understand a remark of \newline
J. Stasheff  on open string theory and higher associative algebras.
\end{quote}
This led to his \cite{kont:trium}.  ``There's an operad for that''; see \cite{ward:cyclic}.
\subsection{Open-Closed Homotopy Algebra and string field theory}\label{ocha}
\h Having considered both \AA- and \LL-algebras, plain and fancy, we come to the combination known as \emph{OCHA} for Open-Closed Homotopy Algebra \cite{hiro-jds:inspired,hiro-jds:mp}. Inspired by open-closed string field theories \cite{zwiebach:mixed}, OCHAs involve an \LL-algebra acting by derivations (up to strong homotopy) on an \AA-algebra but having an additional piece of structure corresponding to a closed string opening to an open string. The relevant operad (rather a colored operad with two colors-one for open and one for closed)
is known as the ``Swiss cheese operad''  (see graphics in \cite{sasha:cheese} for explanation of the name). The details are quite complicated in the original papers, but, just as other
 ``$\infty$'' algebras can be characterized by a single coderivation on an appropriate dgc coalgebra, the same has been achieved for OCHAs by Hoefel \cite{hoefel:coalg}.   A small example appears in \cite{lada-kade:ocha}.
\p
Carqueville has called to my attention the rich class of examples provided by Landau-Ginzburg models \cite{carqueville&kay:invitation,carqueville&kay:otst},
In addition, string theory and string field theory have inspired both string topology, initiated by Chas and Sullivan \cite{chas-sullivan} and a further variety of $\infty$-algebras. 

%\subsection{Feynman ?algebras?}
\subsection{Scattering  amplitudes}\label{scat}
\h Scattering amplitudes in gauge theories are important `observables'. Arkani-Hamed and his colleagues 
\cite{arkani:amplituhedron,arkani:positiveamplituhedron,arkani:scattering} were led to yet another polyhdedron, the \emph{amplituhedron},
which is a generalization of the `positive Grassmannian'.
As we have seen, trees are often a starting point leading to more general structures and tree scattering amplitudes are a good place to start. In \cite{arkani:scattering},
they present a ``novel construction of the associahedron in kinematic space''.

\section{Now and future}
\h That brings us only somewhat  up to date; there is much work in progress ``as we come on the air'' even in this small part of the space of higher structures.
%, espeically in realtion to phyiscs. 
\p
 I find in Manin's \emph{Mathematics, Art, Civilization} \cite{manin:mac}:
\begin{quote}
With the advent of polycategories, enriched categories, \AA-categories, and similar structures, we are beginning to speak a language\dots.
\end{quote}
Now I find this delightfully ironic since, when I first submitted my theses for publication in AJM, they were  deemed too narrow and essentially of no relation to other parts of math!
\p
Perhaps Heraclitus was right: All is flux, nothing stays still. : -)

\end{document}